
\input amstex
\loadbold
\magnification=1100
\def\ra{\rightarrow}
\def\Tor{\operatorname{Tor}}

\def\Ext{\operatorname{Ext}}

\def\dim{\operatorname{dim}}

\def\embdim{\operatorname{embdim}}
\def\depth{\operatorname{depth}}

\def\H{\operatorname{H}}

\define\Qed{\hfill\hbox{\qed}}

\def\m{\operatorname{\frak m}}

\def\C{\operatorname{\bold C}}

\def\F{\operatorname{\bold F}}

\def\G{\operatorname{\bold G}}

\def\C{\operatorname{\bold C}}

\def\Hom{\operatorname{Hom}}
\def\length{\operatorname{length}}

\def\ker{\operatorname{ker}}
\def\tgap{\operatorname{Tor-gap}}
\def\egap{\operatorname{Ext-gap}}
\def\coker{\operatorname{coker}}
\def\image{\operatorname{image}}
\def\ann{\operatorname{ann}}
\def\pd{\operatorname{pd}}

\documentstyle{amsppt}
\pageheight{7 in}
\pagewidth{5 in}
\vcorrection{1 in}
\hcorrection{.5 in}

\rightheadtext{Symmetry in Vanishing Ext}
\topmatter
\title Symmetry in the vanishing of Ext over Gorenstein rings
\endtitle

\author Craig Huneke and David A. Jorgensen\endauthor
\date September 17, 2001 \enddate
\address
Department of Mathematics, University of Kansas,
Lawrence, KS 66045
\endaddress
\email
huneke\@math.ukans.edu
\endemail
\address
Department of Mathematics, University of Texas at Arlington,
Arlington, TX 76019
\endaddress
\email
djorgens\@math.uta.edu
\endemail
\thanks
The first author was partially supported by the NSF and the second
 author was partially supported by the NSA.  This work was
done while the second author was visiting Kansas University.  He
thanks KU for their generous support.
\endthanks

\abstract We investigate 
symmetry in the vanishing of Ext for 
finitely generated modules over local Gorenstein rings.  In particular, 
we define a class of local Gorenstein rings, which we call AB rings, 
and show that for finitely generated modules $M$ and $N$ over an 
AB ring
$R$, $\Ext^i_R(M,N)=0$ for all $i\gg 0$ if and only if 
$\Ext^i_R(N,M)=0$ for all $i\gg 0$.
\endabstract
\endtopmatter

\document

\heading Introduction \endheading

Let $R$ be a local Gorenstein ring   
and let $M$ and $N$ denote
finitely generated $R$-modules.  This paper is concerned with the
relation between the vanishing of all 
higher $\Ext_R(M,N)$ and the vanishing of all higher $\Ext_R(N,M)$.
As a means of investigation we concentrate on the more natural duality 
between the 
vanishing of all higher $\Ext_R(M,N)$ and the vanishing of all higher 
$\Tor$ modules where either $M$ or $N$ is replaced by its dual
$M^*(:=\Hom_R(M,R))$ or $N^*$.

Our interest in this topic came about in part from the following
striking result proved recently by Avramov and Buchweitz \cite{AvBu, Thm. III}.
Suppose $M$ and $N$ are finitely generated
 modules over a complete intersection $R$.
Then the following are equivalent:
\roster
\item $\Tor_i^R(M,N)=0$ for all $i\gg 0$
\item $\Ext^i_R(M,N)=0$ for all $i\gg 0$
\item $\Ext^i_R(N,M)=0$ for all $i\gg 0.$
\endroster
Their proof relies heavily on the use of
certain affine algebraic sets associated to $M$ and $N$, 
called support varieties. In their paper \cite{AvBu}, Avramov and
Buchweitz raise the question of what class of rings satisfy these equivalences
for all finitely generated modules $M$ and $N$. They point out this class
lies somewhere between complete intersections and local Gorenstein rings, but
mention that they do not know whether this class is equal to either
complete intersections or Gorenstein. In this paper we introduce a class of
local Gorenstein rings, which we call \it AB rings \rm, and prove that
AB rings satisfy the property that for all finitely generated modules
$M$ and $N$, the following are equivalent (Theorem 4.1):

\roster  
\item $\Ext^i_R(M,N)=0$ for all $i\gg 0$
\item $\Ext^i_R(N,M)=0$ for all $i\gg 0.$
\endroster

Regular local rings are AB rings, and if $R$ is an AB ring, then 
$R/(x_1,...,x_c)$
is also an AB ring whenever $x_1,...,x_c$ is a regular sequence 
(see Proposition 3.2). This implies
complete intersections are AB rings. Even when restricted to the case of a
complete intersection, our proof of the above equivalence avoids 
the use of support
varieties, and in some ways is more direct than the methods of \cite{AvBu}.
We also prove that local Gorenstein
rings of minimal possible multiplicity are AB rings, for the strong
reason that over such rings (except when the embedding
dimension is 2) all large Ext$_R(M,N)$ vanish if and only if either
$M$ or $N$ has finite projective dimension. See Theorem 3.5 for a precise 
statement. These rings are not complete intersections in general, 
so that in particular the class of AB rings is strictly larger 
than that of complete intersections.

An AB ring $R$ is a local Gorenstein ring defined by the property
that there is a constant $C$,
depending only on the ring, such that if Ext$^i_R(M,N) = 0$ for all
$i \gg 0$, then Ext$^i_R(M,N) = 0$ for all $i>C$. As far as we know
\it every \rm Gorenstein ring is an AB ring; we have been unable to
find an example which is not.  The name `AB' stands for both
Auslander-Bridger and Avramov-Buchweitz.

The paper is organized as follows. In Section 1 we give some preliminary
and straightforward results concerning the relationship
of Ext and Tor. In Section 2 we prove a basic result concerning what
holds over an arbitrary local Gorenstein ring. Specifically, if $M$ and
$N$ are finitely generated maximal Cohen-Macaulay modules over a 
local Gorenstein ring $R$, then the following are equivalent:

\roster
\item $\text{Tor}^R_i(M,N)=0$ \text{ for all $i \gg 0$ },
\item $\text{Ext}^i_R(M,N^*)=0$ \text{ for all $i \gg 0$}, and
\item $\text{Ext}^i_R(N,M^*)=0$ \text{ for all $i \gg 0$}.
\endroster

Section 3 introduces AB rings, details their basic
properties, and gives the main examples. In Section 4 we prove the
main theorem of  symmetry in the vanishing of Ext over AB rings. 
Section 5 contains some independent observations concerning
what the vanishing of Ext means. In particular we relate the vanishing
of $d$ consecutive Ext modules ($d$ being the dimension of the ring)
to the Cohen-Macaulay property of a related tensor product.
We include some questions in a final section.
\bigskip

\heading 1. Preliminaries \endheading

In this section we set notation and discuss some basic facts
which will be used throughout the paper.

Unless otherwise stated, we will assume $R$ to be a local
Gorenstein ring.  
Also, $M$ and $N$ will denote finitely generated $R$-modules.
For an $R$-module $M$ we let $M^*$ denote its dual $\Hom_R(M,R)$.
If $M$ is maximal Cohen-Macaulay then it is also reflexive,
meaning $M^{**}\simeq M$ (assuming $R$ is Gorenstein).

By a {\it complete intersection\/} we mean a local ring whose
completion with respect to the maximal ideal is the quotient 
of a regular local ring by a regular sequence.

For a local ring $R$, we let $\embdim R$ denote its 
{\it embedding dimension\/}.

\subhead Syzygies and Conversions for Ext and Tor \endsubhead

Suppose $M$ is an $R$-module.  Then for $i\geq 0$ 
we let $M_i$ denote $\image f_i$, where $f_i$ is the $i$th
differential in a 
minimal free resolution
$$
\F: \qquad \cdots \to F_2 @>f_2>> F_1 @>f_1>> F_0 @>f_0>> M \to 0
$$
of $M$.  These $M_i$ are the {\it non-negative syzygies\/} of $M$. 
They are unique up to isomorphism, or if one considers a non-minimal
resolution any two are stably isomorphic.

Now suppose that $M$ is a maximal Cohen-Macaulay $R$-module.  Let
$$
\G: \qquad \cdots \to G_2 @>g_2>> G_1 @>g_1>> G_0 @>>> M^* \to 0
$$
a minimal free resolution of its dual $M^*$.  
Since $M^*$ is maximal Cohen-Macaulay,
the dual sequence
$$
\G^*: \qquad  0 \to M^{**} \to G_0^* @>g^*_1>> G_1^* @>g^*_2>> G_2^* 
@>>> \cdots
$$
is exact.  Using the fact that $M$ is reflexive, 
we can splice $\F$ and
$\G^*$ together, getting the doubly infinite long exact sequence
$$
\C(M): \qquad \cdots \to\underset 2\to F_2 @>f_2>> \underset 1\to F_1 
@>f_1>>\underset 0\to  F_0 \to 
\underset {-1}\to G_0^* @>g_1^*>> \underset {-2}\to G_1^* 
@>g_2^*>> \underset {-3}\to G_2^* \to \cdots.
$$
(Note the degree convention.) For $i \leq -1$ we set
$M_i:=\image(g_{-i}^*)$ 
These are the {\it negative syzygies\/} of $M$.  They are unique 
up to isomorphism.  Note that $M_i$ is 
again maximal Cohen-Macaulay for all $i$ when $M$ is.  

We now list some properties of the long exact sequences $\C(M)$.

\medskip

\proclaim{1.1 Lemma} Let $M$ be a finitely generated maximal 
Cohen-Macaulay $R$-module,
and let $N$ be a finitely generated $R$-module.
\roster
\item $\C(M)^* \simeq \C(M^*)$
\item $(M_i)^* \simeq (M^*)_{-i}$ for all $i$.
\item $\Tor^R_i(M,N) \simeq \H_i(\C(M) \otimes N)$ for $i \geq 1$.
\item $\Ext^i_R(M,N) \simeq \H_{-i-1}(\C(M)^* \otimes N)$ for $i \geq 1$.
\item for fixed $t\geq 3$ and for $1\leq i \leq t-2$ we have
\medskip
\noindent$(i)\qquad\Ext^i_R(M_{-t},N)\simeq\Tor^R_{t-i-1}(M^*,N)$

\noindent and

\noindent$(ii)\qquad\Tor_i^R(M_{-t},N)\simeq\Ext_R^{t-i-1}(M^*,N).$
\endroster
\endproclaim

\demo{Proof} Conditions (1)--(3) are straightforward. 
Condition (5) follows easily from (1)--(4), so only Condition (4)
needs some explanation.  The critical fact we need to show is that 
for any complex of free 
$R$-modules $\F$, $\Hom_R(\F,N)$ and
$\Hom_R(\F,R)\otimes_R N$ are isomorphic as complexes:
write 
$$
\F:\qquad\cdots\to F_{i+1} @>f_{i+1}>> F_i @>f_i>> F_{i-1} \to \cdots,
$$
where the $F_i$ are free $R$-modules.  The natural
maps $h_i:\Hom_R(F_i,R)\otimes_R N \to \Hom_R(F_i,N)$ given by
$f\otimes n \mapsto \{a\mapsto f(a)n\}$ are isomorphisms
since $F_i$ is free.  It is easy to check that the diagram
$$
\CD
\Hom_R(F_i,R)\otimes_R N @>f_{i+1}^*\otimes N>> 
\Hom_R(F_{i+1},R)\otimes_R N\\
@VVh_iV @VVh_{i+1}V\\
\Hom_R(F_i,N) @>\Hom(f_{i+1},N)>> \Hom_R(F_{i+1},N)
\endCD
$$
is commutative, and this establishes our fact.
\Qed
\enddemo

\medskip

Suppose $M$ and $N$ are $R$-modules with $M$ maximal Cohen-Macaulay.  
Then $\Ext^i(M,R)=0$ for all $i\geq 1$, and so by shifting along the 
short exact sequences $0\to N_{n+1}\to G_n \to N_n\to 0$
(with $G_n$ free) we obtain isomorphisms
$$
\Ext^i_R(M,N)\simeq\Ext^{i+n}_R(M,N_n)\tag \bf 1.2 \rm
$$
for $i\geq 1$ and $n\geq 0$.

\subhead The Change of Rings Long Exact Sequences of Ext and Tor \endsubhead

Suppose that $S$ is a commutative ring, $x$ is a non-zerodivisor of $S$
and $R:=S/(x)$.  Let $M$ and $N$ be $R$-modules.  Then we have the 
{\it change of rings long exact sequence of $\Ext$\/} \cite{Ro, 11.65}
$$
\matrix
\vdots     & {}    & \vdots     & {}   & \vdots     & {}    & {}\\
\Ext^1_R(M,N)  & @<<<  & \Ext^2_S(M,N)  & @<<< & \Ext^2_R(M,N)  & 
@<<<  & {}\\
\vspace{10pt}
\Ext^0_R(M,N)  & @<<<  & \Ext^1_S(M,N)  & @<<< & \Ext^1_R(M,N) & @<<<  & 0,
\endmatrix\tag 1.3
$$
and the {\it change of rings long exact sequence of $\Tor$\/} \cite{Ro. 11.64}
$$
\matrix
\vdots     & {}    & \vdots     & {}   & \vdots     & {}    & {}\\
\Tor_1^R(M,N)  & @>>>  & \Tor_2^S(M,N)  & @>>> & \Tor_2^R(M,N)  & 
@>>>  & {}\\
\vspace{10pt}
\Tor_0^R(M,N)  & @>>>  & \Tor_1^S(M,N)  & @>>> & \Tor_1^R(M,N) & @>>>  & 0.
\endmatrix\tag 1.4
$$

\bigskip

\heading 2. Vanishing of Ext and Tor over Arbitrary Local Gorenstein Rings
\endheading

\medskip

In this section we prove what type of duality between the vanishing of
Ext and Tor holds over arbitrary local Gorenstein rings. It is possible
an even stronger result is true, as we discuss in Section 4, but
the main result of this section is what is true `on the surface'. In
particular, Theorem 2.1 states that one can flip the arguments in 
vanishing Ext modules `up to duals'.

\medskip

\proclaim{2.1 Theorem} Let $R$ be a local Gorenstein ring, and let $M$ 
and $N$ be finitely generated maximal 
Cohen-Macaulay $R$-modules. Then the following are equivalent: 
{\roster
\item $\Tor^R_i(M,N)=0$ \text{ for all $i \gg 0$},
\item $\Ext^i_R(M,N^*)=0$ \text{ for all $i \gg 0$}, and
\item $\Ext^i_R(N,M^*)=0$ \text{ for all $i \gg 0$}.
\endroster}
\endproclaim

\demo{Proof} Suppose we have shown that (1) and (2) are equivalent. 
By replacing
(1) by the equivalent condition that Tor$^R_i(N,M) = 0$ for all
$i \gg 0$, we see then that (1) is equivalent to (3). Hence
it suffices to prove (1) and (2) are equivalent, and for this
we only need to assume that $N$ is maximal Cohen-Macaulay.

We induce upon the dimension of $R$, say $d$.  If $d = 0$, then 
Ext$^i_R(M,N^*)$ is the Matlis dual of Tor$^R_i(M,N)$,  so the result
is immediate in this case. 

Now assume that $d > 0$. Choose $x\in R$
a non-zerodivisor on $M_d$, $N$, $N^*$, and $R$.  We have
$$
\alignat 2
{}&\quad\Tor^R_i(M,N)=0 &&\qquad\text{for all }i\gg 0;\\
\Longleftrightarrow&\quad\Tor^R_i(M_d,N)=0&&\qquad\text{for all }i\gg 0;\\
\intertext{and from the long exact sequence of Tor coming from the
short exact sequence $0\to N @>x>> N\to N/xN\to 0$ and Nakayama's lemma,}
\Longleftrightarrow&\quad\Tor^R_i(M_d,N/xN)=0&&\qquad\text{for all }i\gg 0;\\
\intertext{and by the standard isomorphisms $\Tor^R_i(M_d,N/xN)\simeq
\Tor^{R/(x)}_i(M_d/xM_d,N/xN)$,}
\Longleftrightarrow&\quad\Tor^{R/(x)}_i(M_d/xM_d,N/xN)=0
&&\qquad\text{for all }i\gg 0;\\
\intertext{by the induction hypothesis,}
\Longleftrightarrow&\quad\Ext_{R/(x)}^i(M_d/xM_d,(N/xN)^*)=0
&&\qquad\text{for all }i\gg 0;\\
\intertext{and now since $N$ is maximal Cohen-Macaulay and $R$ 
is Gorenstein, $N^*/xN^*\simeq (N/xN)^*$ (where the second module is 
$\Hom_{R/(x)}(N/xN,R/(x))$), thus,}
\Longleftrightarrow&\quad\Ext_{R/(x)}^i(M_d/xM_d,N^*/xN^*)=0
&&\qquad\text{for all }i\gg 0;\\
\intertext{by the isomorphisms \cite{R} $\Ext^i_R(M_d,N^*/xN^*)\simeq
\Ext^i_{R/(x)}(M_d/xM_d,N^*/xN^*)$,}
\Longleftrightarrow&\quad\Ext_R^i(M_d,N^*/xN^*)=0
&&\qquad\text{for all }i\gg 0;\\
\intertext{and now from the long exact sequence of Ext coming
from the short exact sequence $0\to N^*@>x>> N^*\to N^*/xN^*\to 0$
and Nakayama's lemma,}
\Longleftrightarrow&\quad\Ext_R^i(M_d,N^*)=0
&&\qquad\text{for all }i\gg 0;\\
\Longleftrightarrow&\quad\Ext_R^i(M,N^*)=0
&&\qquad\text{for all }i\gg 0.\\
\endalignat
$$ 
\Qed
\enddemo

\medskip

\remark{{\rm\bf 2.2} Remark} Suppose $M$ and $N$ are maximal Cohen-Macaulay
modules over the local Gorenstein ring $R$.  Then
$$
\Ext^i_R(M,N)\simeq\Ext^i_R(N^*,M^*).
$$

This isomorphism can be seen as follows:  suppose that $i=1$.
As $M$ and $N$ are maximal Cohen-Macaulay, they are reflexive, so
short exact sequences
$0\ra N\ra T\ra M\ra 0$ dualize to short exact sequences
$0\ra M^*\ra T^*\ra N^*\ra 0$ and vice-versa. The Yoneda definition of
$\Ext^1$ then gives the isomorphism.

For $i>1$ we have 
$$
\alignat 2
\Ext^i_R(M,N)&\simeq\Ext^1_R(M_{i-1},N)&&{}\\
{}&\simeq\Ext^1_R(N^*,M_{i-1}^*)&&\quad\text{by the $i=1$ case}\\
{}&\simeq\Ext^1_R(N^*,(M^*)_{-i+1})&&\quad\text{by (2) of Lemma 1.1}\\
{}&\simeq\Ext^i_R(N^*,M^*)&&\quad\text{by (1.2)}.
\endalignat
$$
\Qed
\endremark

Presumably, one can directly prove this remark using the Yoneda definition
of Ext and the fact that both $M$ and $N$ are maximal Cohen-Macaulay.

\medskip
 
Below is an example showing that the hypothesis that $N$ is
maximal Cohen-Macaulay in the equivalence of (1) and (2) in
Theorem 2.1 cannot be dropped.

\medskip

\example{2.3 Example} Let $R$ be the 3-dimensional hypersurface 
$k[[W,X,Y,Z]]/(WX-YZ)$, and set $M:=k$ and
$N:=\coker\left(\smallmatrix w \\ x \\ y \\ z \endsmallmatrix\right).$
Then $\pd_RN=1$ (but $\pd_RN^*=\infty$), so we have
$\Tor^R_i(M,N)=0$ for all $i\gg 0$.  However, $\Ext^4_R(M,N^*)\neq 0$.
By what is shown in the next section, $R$ is an AB ring. 
If it were the case that $\Ext^i_R(M,N^*)$ is zero for all
$i \gg 0$ then, as $R$ is AB, Proposition 3.1 shows then
that $\Ext^i_R(M,N^*) = 0$ for all $i> \dim R$, which would be
a contradiction.
\Qed\endexample 

\bigskip

\heading 3. AB Rings \endheading
 
\medskip

Let $R$ be a commutative ring.  We define the 
Ext-{\it index\/} of $R$ to be
$$
\sup\{n\vert\Ext^i_R(M,N)=0\text{ for all }i>n\text{ and }
\Ext^n_R(M,N)\neq 0\},
$$
where the sup is taken over all pairs of finitely generated 
$R$-modules $(M,N)$ with $\Ext^i_R(M,N)=0$ for all $i\gg 0$.
 
\definition{3.0 Definition}
If $R$ is a local Gorenstein ring of finite Ext-index, we say that
$R$ is an {\it AB ring\/}.
\enddefinition

We will prove that all complete intersections
are AB rings. More generally, it is obvious that $R$ is an AB ring if 
$\hat R$ is (where $\hat R$ is the completion of $R$), and we show
(3.2) that if $R$ is an AB ring and $x_1,...,x_c$ is
a regular sequence, then $R/(x_1,...,x_c)$ is also an AB ring. The class
of AB rings also includes local 
Gorenstein rings of `minimal' multiplicity
embdim$(R) - \text{dim}(R) + 2$ (see 3.5).

\medskip

\proclaim{3.1 Proposition} Suppose that $R$ is an AB ring. 
Then the $\Ext$-index of $R$ equals $\dim R$.
\endproclaim

\demo{Proof} Let $n$ denote the Ext-index of $R$ and $d$ the dimension 
of $R$.

Let $x_1,\dots,x_d$ be a maximal $R$-sequence.  Set $M:=R/(x_1,\dots,x_d)$
and $N=k$, the residue field of $R$.  Then $\Ext^i_R(M,N)=0$ for $i>d$
and $\Ext^d_R(M,N)\simeq k\neq 0$.  Hence $n\geq d$.

Suppose that $n>d$.  There exists a pair of finitely generated 
$R$-modules $(M,N)$
such that $\Ext^i_R(M,N)=0$ for $i>n$ and $\Ext^n_R(M,N)\neq 0$.  We 
have the isomorphisms 
$\Ext^{i+1}_R((M_d)_{-d-1},N)\simeq\Ext^{i-d}_R(M_d,N)
\simeq\Ext^i_R(M,N)$ for $i>d$.
Hence $\Ext^i((M_d)_{-d-1},N)=0$ for $i>n+1$ and 
$\Ext^{n+1}((M_d)_{-d-1},N)\neq 0$, which contradicts the definition
of $n$.  Therefore $n=d$.
\Qed
\enddemo

We of course have the dual notion of Tor-{\it index\/}.  If $R$
is local Gorenstein with finite Tor-index, then it is also equal to
$\dim R$, by an argument analogous to that of 3.1.

Another related property of rings we are interested in is the following.
We say that $\Ext_R^*(M,N)$ has a {\it gap of length
$t$} if for some $n\geq 0$, $\Ext_R^i(M,N)=0$ for $n+1\leq i\leq n+t$, but
$\Ext_R^n(M,N)$ and $\Ext_R^{n+t+1}(M,N)$ 
are both nonzero.  We have the analogous notion of gap for $\Tor_*^R(M,N)$.
(We allow gaps of length $0$.)
We set
$$
\egap(R):=\sup\left\{t\in\Bbb N\text{ }\bigg\vert\matrix 
\Ext^*_R(M,N) \text{ has a gap of length $t$ for } \\ 
\text{ finite $R$-modules $M$ and $N$}
\endmatrix\right\},
$$
and
$$
\tgap(R):=\sup\left\{t\in\Bbb N\text{ }\bigg\vert\matrix 
\Tor^R_*(M,N) \text{ has a gap 
of length $t$ for } \\ \text{ finite $R$-modules $M$ and $N$ } 
\endmatrix\right\}.
$$

We say that $R$ is Ext-{\it bounded\/} if it has finite Ext-gap.
Similarly, we say $R$ is Tor-{\it bounded\/} if it has finite
Tor-gap.

We list some elementary properties involving finite Ext-index, Tor-index,
Ext-boundedness and Tor-boundedness for local Gorenstein rings.

\proclaim{3.2 Proposition}  Let $x$ be a non-zerodivisor of the 
$d$-dimensional local Gorenstein ring $R$. Then 
\roster
\item $R$ is an AB ring if and only if $R/(x)$ is an AB ring.
\item $R$ has finite Tor-index if and only if $R/(x)$ does.
\item $R$ is Ext-bounded if and only if $R/(x)$ is.
\item $R$ is Tor-bounded if and only if $R/(x)$ is.
 
\endroster
\endproclaim

\demo{Proof} (1).
Suppose that $R$ is an AB ring.  
Let $M$ and $N$ be finitely generated
$R/(x)$-modules such that
$\Ext_{R/(x)}^i(M,N)=0$  for all $i\gg 0$.
By the change of rings long exact sequence of Ext (1.3)
we conclude that $\Ext_R^i(M,N)=0$ for all $i\gg 0$, and so
$\Ext_R^i(M,N)=0$ for all $i>d$.  Looking at (1.3) again, 
we see that $\Ext_{R/(x)}^i(M,N)\simeq\Ext_{R/(x)}^{i+2}(M,N)$ for
$i>d-1$.  But as $\Ext_{R/(x)}^i(M,N)=0$ for all $i\gg 0$, we have
$\Ext_{R/(x)}^i(M,N)=0$ for all $i>d-1$.  Hence $R/(x)$
is an AB ring.

Now suppose that $R/(x)$ is an AB ring, and  
let $M$ and $N$ be finitely generated $R$-modules
such that $\Ext^i_R(M,N)=0$ for all $i\gg 0$.
We have the isomorphisms 
$$
\Ext^i_R(M,N)\simeq\Ext^{i-d}_R(M_d,N)
\simeq\Ext^i_R(M_d,N_d)
$$
which are valid for $i>d$, the second one being that of (1.2).
The short exact sequence 
$0 \to N_d @>x>> N_d \to N_d/xN_d \to 0$ gives rise
to the the long exact sequence of Ext
$$
\cdots \to \Ext^i_R(M_d,N_d) @>x>>
\Ext^i_R(M_d,N_d) \to \Ext^i_R(M_d,N_d/xN_d) \to \cdots.\tag 3.2.1
$$
Since $\Ext^i_R(M_d,N_d)=0$ for all $i\gg 0$, we see that  
$\Ext^i_R(M_d,N_d/xN_d)=0$ for all $i \gg 0$.  We have the 
isomorphisms \cite{R}
$$
\Ext^i_{R/(x)}(M_d/xM_d,N_d/xN_d)\cong\Ext^i_R(M_d,N_d/xN_d)\tag 3.2.2
$$
for all $i\geq 0$.  Hence $\Ext^i_{R/(x)}(M_d/xM_d,N_d/xN_d)=0$
for all $i\gg 0$, which means that 
$\Ext^i_{R/(x)}(M_d/xM_d,N_d/xN_d)=0$ for all $i>d-1$, since
$R/(x)$ is an AB ring.
Therefore $\Ext^i_R(M_d,N_d/xN_d)=0$
for all $i>d-1$.  By (3.2.1) and Nakayama's Lemma, we conclude
that $\Ext^i_R(M_d,N_d)=0$ for all $i>d-1$, and so $\Ext^i_R(M,N)=0$
for all $i>d$.  Therefore $R$ is an AB ring.

The proof of (2) is exactly analogous to the proof of (1), using
(1.4) and a long exact sequence of Tor this time.

(3).  Assume that $e:=\egap(R)<\infty$. Let $M$ and $N$ be 
finitely generated $R/(x)$-modules such that $\Ext_{R/(x)}^i(M,N)=0$ for 
$n\leq i\leq n+e+1$, some $n\geq 1$.  
The change of rings long exact sequence of Ext (1.3)
shows that $\Ext_R^i(M,N)=0$ for $n+1\leq i\leq n+e+1$.  Since
$\egap(R)=e$ we have $\Ext_R^i(M,N)=0$ for all $i>n$.  Another look
at (1.3) shows that 
$\Ext_{R/(x)}^i(M,N)\simeq\Ext_{R/(x)}^{i+2}(M,N)$ for all $i>n-1$.
Since $\Ext_{R/(x)}^i(M,N)=0$ for $i=n,n+1$, we see then that
$\Ext_{R/(x)}^i(M,N)=0$ for $i>n-1$
Hence $\egap(R/(x))\leq e+1$.

Now assume that $e:=\egap(R/(x))<\infty$.  Suppose that $M$ and $N$
are finitely generated $R$-modules with 
$\Ext_R^i(M,N)=0$ for $n\leq i\leq n+d+e+1$, some $n\geq 1$.  
We have $\Ext_R^i(M_d,N_d)\simeq\Ext_R^i(M,N)=0$ for
$n+d\leq i\leq n+d+e+1$.
Therefore, from (3.2.1), we get
$\Ext^i_R(M_d,N_d/xN_d)=0$ for $n+d\leq i\leq n+d+e$.  Equivalently,
$\Ext_{R/(x)}^i(M_d/xM_d,N_d/xN_d)=0$ for $n+d\leq i\leq n+d+e$.
Since $\egap(R/(x))=e$, $\Ext_{R/(x)}^i(M_d/xM_d,N_d/xN_d)=0$
for all $i\geq n+d$, which implies, by (3.2.1)
and Nakayama's lemma, $\Ext_R^i(M_d,N_d)=0$ for all $i\geq n+d$,
which means $\Ext_R^i(M,N)=0$ for all $i\geq n$. Therefore
$\egap(R)\leq d+e+1$.

The proof of (4) is similar to the proof of (3).
\Qed\enddemo

\medskip

\proclaim{3.3 Theorem} Assume that $R$ is a local Gorenstein ring.
Then 
\roster
\item $R$ is an AB ring if and only if it has finite Tor-index;
\item $R$ is Ext-bounded if and only if it is Tor-bounded;
\item if $R$ is Ext-bounded, then it is an AB ring.
\endroster
\endproclaim

\demo{Proof}  We first prove (1). Choose a maximal regular sequence
in $R$ and let $I$ be the ideal generated by this sequence. Proposition
3.2 states that $R$ is an AB ring if and only if $R/I$ is an AB ring,
and $R$ has finite Tor-index if and only if $R/I$ has finite Tor-index.
Hence it suffices to prove (1) in case $R$ is $0$-dimensional. In this case
$\Ext^i_R(M,N)$ is the Matlis dual of $\Tor^R_i(M,N^*)$ so that the vanishing
of one implies the vanishing of the other. This proves (1).

Statement (2) is proved in a similar manner, using Proposition 3.2.

We prove (3).  Assume $R$ is Ext-bounded.
We prove that $R$ has finite Tor-index.  
Let $d$ denote the dimension
of $R$ and $e:=\egap(R)$, and suppose that
for finite $R$-modules $M$ and $N$, $\Tor^R_i(M,N)=0$ for all $i\gg 0$. 
Let $b=d-\depth M$ so that $M_b$ is maximal Cohen-Macaulay.  Choose 
$n$ largest such that $\Tor_n^R(M_b,N)\neq 0$. Using $t=e+n+3$. 
in $(5)(i)$ of Lemma 1.1, we have 
$\Ext^i_R((M_b^*)_{-e-n-3},N)\simeq\Tor^R_{e+n+2-i}(M_b,N)=0$
for $1\leq i\leq e+1.$  Hence we have a gap of zero Ext
larger than $e$.  Therefore $\Ext_R^i((M_b^*)_{-e-n-3},N)=0$ for all 
$i\geq 1$, which forces $n=0$.  
Thus $\Tor_i^R(M,N)=0$ for all $i>d$.
\Qed
\enddemo

\medskip
The following Corollary is an almost immediate consequence of Proposition
3.3, as regular local rings are clearly Ext-bounded.

\proclaim{3.4 Corollary} Let $R$ be a local Gorenstein ring. If $R$ is
a complete intersection, then $R$ is Ext-bounded. In particular,
$R$ is an AB ring.
\endproclaim

\demo{Proof} Since $R \hookrightarrow \hat R$ is a faithfully
flat extension, $\Ext^i_R(M,N)=0$ if and only if 
$\Ext^i_{\hat R}(\hat M,\hat N)=0$ and so $R$ is an AB ring if
$\hat R$ is.  Therefore we may 
without loss of generality assume that
$R=S/(x_1,\dots,x_c)$ where $S$ is a regular local ring and 
$x_1,\dots,x_c$ is an $S$-regular sequence.  By Proposition 3.2
it suffices to prove that $S$ is Ext-bounded.  But this is trivial
as every finitely generated module over $S$ has projective 
dimension $\leq\dim S$,
so that $\Ext^i_R(M,N)=0$ for $i>\dim S$ and Ext-gaps can occur of 
length no longer than $\dim S-2$.
\Qed
\enddemo

\medskip

It could be that {\it all \/} local Gorenstein rings are AB rings; 
we have no
counterexample. The class of AB rings is strictly bigger than
the class of complete intersections as the next Theorem proves,
albeit for rather strong reasons.

\medskip

\proclaim{3.5 Theorem} Let $(R,\m,k)$ be a local Gorenstein ring with
multiplicity equal to $\embdim R-\dim R+2$.  Assume that $\embdim R>2$
(so that $R$ is not a complete intersection).  Then for finitely generated 
$R$-modules
$M$ and $N$, $\Ext^i_R(M,N)=0$ for all $i\gg 0$ if and only if
either $M$ or $N$ has finite projective dimension. In particular,
$R$ is an AB ring.
\endproclaim

\demo{Proof} We induce on $d:=\dim R$.

$d=0$.  In this case by duality we have that $\Tor_i^R(M,N^*)=0$
for all $i\gg 0$.  Replace $M$ by a high enough syzygy such that
$\Tor_3^R(M,N^*)=\Tor_4^R(M,N^*)=\Tor_5^R(M,N^*)=0$.  Now the 
following lemma, 3.6, says that either $M$ or $N$ is free.

$d>0$.  By replacing $M$ and $N$ by syzygies we can assume they 
are both maximal Cohen-Macaulay.  Let $x$ be a non-zerodivisor
of $R$.  Once again we use the fact that $\Ext^i_R(M,N)=0$
for all $i\gg 0$ if and only if $\Ext^i_{R/(x)}(M/xM,N/xN)=0$
for all $i\gg 0$.  By induction either $M/xM$ or $N/xN$ has finite
projective dimension over $R/(x)$.  But then either $M$ or $N$
has finite projective dimension over $R$.
\Qed
\enddemo

\medskip

\proclaim{3.6 Lemma} Let $(R,\m,k)$ be a $0$-dimensional
local Gorenstein ring with multiplicity $\embdim R+2$.  
Assume $\embdim R>2$ (so that $R$ is not a complete
intersection).  Let
$M$ and $N$ be finitely generated $R$-modules.  Then
$\Tor^3_R(M,N)=\Tor^4_R(M,N)=\Tor^5_R(M,N)=0$
if and only if either $M$ or $N$ is free.
\endproclaim

\demo{Proof} Let $n$ denote the embedding dimension of $R$.
Assume that $M$ is not free.
If $k$ is a summand of $M_i$ for any 
$0\leq i \leq 4$, then we get right away that $N$ is free,
since $\Tor^R_1(M_i,N)=\Tor^R_{i+1}(M,N)=0$ would then imply
$\Tor^R_1(k,N)=0$.
Therefore assume $k$ is not a summand of $M_i$ for all 
$0\leq i \leq 4$.  Replace $M$ by its first syzygy.  Then
as $M\subseteq \m F$, for $F$ a free module, we have $\m^2M=0$
(since $\m^3=0$).
Let $b_i$ denote the $i$th Betti number of $M$ and $s:=\dim_k\m M$.
Then, as in Lescot's paper \cite{L, Lemma 3.3}, $M$ is $3$-exceptional
and $b_1=nb_0-s$, $b_2=b_0(n^2-1)-sn$ and $b_3=b_0(n^3-2n)-s(n^2-1)$.

Now suppose that $N$ is also not free.  Also replace $N$ by its first syzygy,
so that $\m^2 N=0$.  Write $\m N\simeq k^d$.  We have a short
exact sequence $0\to k^d\to N\to k^c\to 0$, where $c$ is the
minimal number of generators of $N$.  
Applying $M\otimes_R\underline{\text{\ \ }}$
to this short exact sequence and using the fact that 
$\Tor^1_R(M,N)=\Tor^2_R(M,N)=\Tor^3_R(M,N)=0$ we get
$cb_2=db_1$ and $cb_3=db_2$.  Letting $\alpha:=d/c$ we can write these
equations as $b_2=\alpha b_1$ and $b_3=\alpha b_2=\alpha^2 b_1$. 
Substituting for the $b_i$ we get
$$
\align
b_0(n^2-1)-sn&=\alpha(nb_0-s)\\
b_0(n^3-2n)-s(n^2-1)&=\alpha^2(nb_0-s).
\endalign
$$
After rearranging we arrive at
$$
\align
b_0(n^2-\alpha n-1)&=s(n-\alpha)\\
b_0(n^3-\alpha^2 n-2n)&=s(n^2-\alpha^2-1).
\endalign
$$
Now cross multiplying, cancelling off the $b_0 s$ terms, and simplifying
we are left with the condition $\alpha^2-n\alpha +1=0$.
This says that $\alpha \in \Bbb Q$ is algebraic over $\Bbb Z$.  Hence
$\alpha$ is an integer, and the only choice is $\alpha=1$.  But then
$n=2$, which is a contradiction.  Hence $N$ must be free.
\Qed
\enddemo

\bigskip

\heading 4. Vanishing of Ext and Tor over AB Rings \endheading

\medskip

In this section we prove that AB rings are a class which gives
the duality of vanishing Ext discussed in the introduction. Our
main theorem states: 

\medskip

\proclaim{4.1 Theorem} Suppose that  $R$ is an AB ring,
and let $M$ and $N$ be finitely generated $R$-modules. Then
$$
\align
{}&\Ext^i_R(M,N)=0 \text{ for all $i\gg 0$ } \text{ if and only if}\\
{}&\Ext^i_R(N,M)=0 \text{ for all $i\gg 0$}.
\endalign
$$ 
\endproclaim

\demo{Proof} First assume the theorem is true if both $M$ and $N$ are
maximal Cohen-Macaulay.  For the general case, take syzygies
$M_m$ and $N_n$ ($m,n\geq 0$) of $M$ and $N$, respectively, which are 
maximal Cohen-Macaulay.  We have $\Ext_R^i(M,N)=0$ for all $i\gg 0$
if and only if $\Ext_R^i(M_m,N)=0$ for all $i\gg 0$ and by (1.2)
this is equivalent to $\Ext^i_R(M_m,N_n)=0$ for all $i \gg 0$. Thus
$\Ext_R^i(M,N)=0$ for all $i\gg 0$ if and only if 
$\Ext^i_R(M_m,N_n)=0$ for all $i\gg 0$, and so the theorem 
holds generally.

Now suppose that $M$ and $N$ are maximal Cohen-Macaulay and
$\Ext_R^i(M,N)=0$ for all $i \gg 0$. Then for all
$t\geq 1$, $\Ext^i_R(M_{-t},N)=0$ for all $i \gg 0$. Since $R$ is
an AB ring, it follows from Proposition 3.1 that for all $t\geq 1$ 
and all $i>d:=\dim (R)$,
$\Ext^i_R(M_{-t},N)=0$. However, 
$\Ext^i_R(M_{-t},N)\simeq\Ext^1_R(M_{i-t-1},N).$
Hence for all $t\geq 1$ and all $i>d$,
$\Ext^1_R(M_{i-t+1},N) = 0.$
By varying $i$ and $t$, we obtain that $\Ext^1_R(M_{-t},N)=0$ for all
$t\geq 1$.
Therefore by $(5)(i)$ of Lemma 1.1,
$\Tor^R_{t-2}(M^*,N) = \Tor^R_{t-2}(N, M^*) =0$ for all $t\geq 3$.
Applying Theorem 2.1 then shows that $\Ext_R^i(N,M)=0$ for all $i\gg 0$.
\Qed
\enddemo

\medskip

As an immediate corollary, we have an analogue of Theorem 2.1

\medskip

\proclaim{4.2 Corollary} Suppose that $R$ is an AB ring,
and let $M$ and $N$ be finitely generated maximal Cohen-Macaulay
$R$-modules.  Then the following are equivalent:
\roster
\item $\Tor_i^R(M,N)=0$ for all $i \gg 0$,
\item $\Ext^i_R(M^*,N)=0$ for all $i\gg 0$, and
\item $\Ext^i_R(N^*,M)=0$ for all $i\gg 0$.
\endroster
\endproclaim

\demo{Proof} Due to the natural symmetry in Tor, it suffices
to prove just the equivalence between (1) and (2), and for this we
only need to assume that $M$ is maximal Cohen-Macaulay. 

We have 
$$
\alignat 2
\Tor^R_i(M,N)=0\text{ for all $i\gg 0$}&\Longleftrightarrow\Tor^R_i(N,M)=0
\text{ for all $i\gg 0$}&&{}\\
{}&\Longleftrightarrow\Ext^i_R(N,M^*)=0
\text{ for all $i\gg 0$}&&\qquad\text{by 2.1}\\
{}&\Longleftrightarrow\Ext^i_R(M^*,N)=0
\text{ for all $i\gg 0$}&&\qquad\text{by 4.1.}
\endalignat
$$
\Qed
\enddemo

\medskip

Our final proposition in this section is an observation that
there are circumstances other than where one module
has finite projective dimension or where the ring is a complete
intersection in which all large Ext modules vanish. 

\medskip

\proclaim{4.3 Proposition-Example} Let $(R,\m_R,k)$ and $(S,\m_S,k)$
be two local Gorenstein rings essentially of finite type over the
same field $k$, and let $M_R$ be a finitely generated $R$-module and
$N_S$ a finitely generated $S$-module. Set $A:=(R\otimes_kS)_P$ where
$P:=(\m_R\otimes_kS+R\otimes_k\m_S)$, $M:=(M_R\otimes_kS)_P$
and $N:=(R\otimes_kN_S)_P$.  Then $A$ is local Gorenstein and
$\Ext^i_A(M,N)=0$ for all $i>\dim A$.
\endproclaim
 
\demo{Proof}  Of course $A$ is Noetherian, being a localization of a
finitely generated $k$-algebra. $A$ is also Gorenstein by applying
\cite{WITO}. 

For the last statement, we induce on $d:=\dim A=\dim R+\dim S$.
Suppose that $A$ has dimension $0$.  In this case duality yields
$\Ext^i_A(M,N)^* \simeq \Tor^A_i(M,N^*)$.  Thus it suffices to prove
$\Tor^A_i(M,N^*)=0$ for all $i>0$.

We first claim that 
$\Hom_{R\otimes_k S}(M_R\otimes_k S,R\otimes_k S)
\simeq\Hom_R(M_R,R)\otimes_k S$.
To see this, note that these modules are naturally isomorphic if $M_R$ is 
free $R$-module.  In general, let $R^m @>\rho>> R^n \to M_R\to 0$
be a presentation of $M_R$ over $R$.  Let $A':=R\otimes_k S$ and
$M':=M_R\otimes_k S$.  This yields a presentation
$(A')^m @>>> (A')^n\to M'\to 0$ of $M'$
over $A'$.  We obtain a commutative diagram
$$
\minCDarrowwidth{20 pt}
\CD
0@>>>\Hom_{A'}(M',A') @>>> \Hom_{A'}((A')^n,A') @>>> \Hom_{A'}((A')^m,A')\\
@. @. @VVV @VVV\\
0@>>>\Hom_R(M_R,R)\otimes_k S @>>> \Hom_R(R^n,R)\otimes_k S
@>>> \Hom_R(R^m,R)\otimes_k S,
\endCD
$$
where the first row is exact
and the vertical arrows are isomorphisms.  To establish the claim we
only need to know that the bottom row is exact, but this follows from
the fact that $0\to\Hom_R(M_R,R)\to\Hom_R(R^n,R)\to
\Hom_R(R^m,R)$ is an exact sequence of $k$-modules and
$S$ is flat as a $k$-module.

Localizing the isomorphism in the claim above at $P$, we see that 
the $A$-module $M^*:=\Hom_A(M,A)$ comes from the $R$-module
$\Hom_R(M_R,R)$.  Similarly $N^*:=\Hom_A(N,A)$ comes from 
the $S$-module $\Hom_S(N_S,S)$.  Hence there is no distinction 
between proving $\Tor^A_i(M,N^*)=0$ for all $i>0$ and proving
$\Tor^A_i(M,N)=0$ for all $i>0$.  We will prove the latter.
 
Let $(\F,f)$ be an $R$-free resolution of $M_R$. 
Then $\F$ is an
exact sequence of $k$ modules, and since $S$ is flat as a $k$-module,
$\F \otimes_k S$ is an exact sequence, of $R\otimes_k S$-modules.  
Thus $(\F\otimes_k S)_P$ is an $A$-free resolution of $M$.  To show
that $\Tor^A_i(M,N)=0$ for all $i>0$ we will simply show that the
complex $(\F\otimes_k S)_P\otimes_A N$ is acyclic (meaning the homology
is zero except in degree zero).

For all $i$ we have a commutative diagram
$$
\CD
(F_i\otimes_kS)\otimes_{R\otimes_kS}(R\otimes_kN_S)
@>(f_i\otimes S)\otimes(R\otimes N_S)>>
(F_{i-1}\otimes_kS)\otimes_{R\otimes_kS}(R\otimes_kN_S)\\
@VV\simeq V @VV\simeq V\\
F_i\otimes_k N_S 
@>f_i\otimes N_S >>
F_{i-1}\otimes_k N_S,
\endCD
$$
where the vertical arrows are the natural isomorphisms.
Hence $(\F\otimes_k S)\otimes_{R\otimes_kS} (R\otimes_k N_S)$ and
$\F\otimes_k N_S$ are isomorphic complexes of ${R\otimes_kS}$-modules.
Since $N_S$ is flat as a $k$-module, the latter is acyclic, and
therefore so is $(\F\otimes_k R_S)\otimes_{R\otimes_kS} (R\otimes_k N_S)$.
Finally, localizing at $P$ we get that $(\F\otimes_kS)_P\otimes_A N$
is acyclic, and this finishes the proof in the $d=0$ case.

Now without loss of generality assume that $\dim R>0$.  
>From the discussion above we know that $M_1\simeq((M_R)_1\otimes_kS)_P$.
Let $x$ be a non-zerodivisor on both $(M_R)_1$ and $R$.
Then $x\otimes 1$ is a non-zerodivisor on $M_1$, $A$ and $N$, and we have
$$
A/(x\otimes 1)\simeq (R/(x)\otimes_k S)_P,
$$
$$
M_1/(x\otimes 1)M_1\simeq ((M_R)_1/x(M_R)_1\otimes_k S)_P
$$ 
and 
$$
N/(x\otimes 1)N\simeq (R/(x)\otimes_k N_S)_P.
$$
Hence by induction we have that 
$$
\Ext^i_{A/(x\otimes 1)}(M_1/(x\otimes 1)M_1,N/(x \otimes 1)N)=0
$$ 
for all $i>d-1$.
Now (3.2.1) and (3.2.2) show that 
$\Ext^i_A(M_1,N)=0$ for all $i>d-1$, which means that 
$\Ext^i_A(M,N)=0$ for all $i>d$.
\Qed
\enddemo

\bigskip
 
\heading 5. What does the vanishing of Ext mean? \endheading

\medskip
 
Many of the results in this section are closely related to the
work of Auslander and Bridger. See \cite{AB}, and the
writeup \cite{M} of the contents of \cite{AB}. However, none of the
results below is explicitly in these works, and we found they gave us
a better understanding of what the vanishing of Ext means.
 
\medskip

\subhead The natural maps $M^*\otimes_R N\to\Hom_R(M,N)$ 
and $M\otimes_R N^*\to\Hom_R(M,N)^*$ \endsubhead

Assume that $M$ is maximal Cohen-Macaulay.  From the short exact sequence
$0\to M_1\to F\to M\to 0$ we get the dual short exact sequence
$0\to M^*\to F^*\to M_1^*\to 0$, and these yield a
commutative diagram
$$
\minCDarrowwidth{20pt}
\CD
{} @. M^*\otimes_R N @>\alpha>> F^*\otimes_R N
@>>> M_1^*\otimes_R N @>>> 0\\
@. @VVf_0V @VVgV @VVf_1V @. \\
0 @>>> \Hom_R(M,N) @>>> \Hom_R(F,N) @>\beta>>
\Hom_R(M_1,N), @. {}
\endCD
$$
where the vertical arrows are the natural maps $M^*\otimes_RN\to\Hom_R(M,N)$
given by $\phi\otimes n \mapsto\{m\mapsto\phi(m)n\}$.  Note that $g$ is an 
isomorphism (since $F$ is free), $\ker\alpha\simeq\Tor^R_1(M_1^*,N)$ and
$\coker\beta\simeq\Ext^1_R(M,N)$.
>From this diagram one easily deduces the following three facts.
\roster
\item $\ker f_1\simeq \coker f_0$,
\item $\Tor^R_1(M_1^*,N)\simeq \ker f_0$, and
\item $\Ext^1_R(M,N)\simeq \coker f_1$.
\endroster
 
Building a diagram as above for each of the exact sequences
$0\to M_{i+1}\to F_i\to M_i\to 0$ and using the corresponding
three facts as above, we see that we have exact
sequences
$$
0\to\Ext^1_R(M_{i-2},N)\to M_i^*\otimes_R N\to\Hom_R(M_i,N)
\to\Ext^1_R(M_{i-1},N)\to 0,
$$
and
$$
0\to\Tor_1^R(M_{i+1}^*,N)\to M_i^*\otimes_R N\to\Hom_R(M_i,N)
\to\Tor^R_1(M_{i+2}^*,N)\to 0.
$$
For $i\geq 2$ the first exact sequence can be written as
$$
0\to\Ext^{i-1}_R(M,N)\to M_i^*\otimes_R N\to\Hom_R(M_i,N)
\to\Ext^i_R(M,N)\to 0.\tag \bf 5.1 \rm
$$
An immediate observation is

\proclaim{5.2 Proposition} Let $R$ be a local Gorenstein ring, and
let $M$ and $N$ be finitely generated $R$-modules with $M$ maximal 
Cohen-Macaulay.  Then
$\Ext^i_R(M,N)=0$ for all $i\gg 0$ if and only if 
the natural maps $M_i^*\otimes_R N\to\Hom_R(M_i,N)$ are isomorphisms
for all $i\gg 0$.\Qed
\endproclaim

Note also that building exact sequences (5.1) for arbitrarily large negative
syzygies of $M$, and
then splicing the resulting exact sequences together, we obtain
a doubly infinite long exact sequence
$$
\cdots\to M^*_i\otimes_R N \to \Hom_R(M_i,N)\to
M^*_{i+i}\otimes_R N \to \Hom_R(M_{i+1},N)\to\cdots. \tag \bf 5.3 \rm
$$
 
\medskip

Now suppose that $N$ is maximal Cohen-Macaulay and that $\Ext^1_R(M,N)=0$.
>From the short exact sequence $0\to M_1 \to F\to M\to 0$ we get the
short exact sequence $0\to\Hom_R(M,N)\to\Hom_R(F,N)\to\Hom_R(M_1,N)\to 0$,
and a commutative diagram
$$
\minCDarrowwidth{20pt}
\CD
{} @. M_1\otimes_R N^* @>\gamma>> F\otimes_R N^*
@>>> M\otimes_R N^* @>>> 0\\
@. @VVh_1V @VVg'V @VVh_0V @. \\
0 @>>> \Hom_R(M_1,N)^* @>>> \Hom_R(F,N)^* @>\delta>>
\Hom_R(M,N)^*, @. {}
\endCD
$$
where the vertical arrows are the natural maps $M\otimes_RN^*\to\Hom_R(M,N)^*$
given by $m\otimes\phi\mapsto\{\psi\mapsto\phi(\psi(m))\}$.   
Note that $g'$ is an 
isomorphism (since $F$ is free), $\ker\gamma\simeq\Tor^R_1(M,N^*)$ and
$\coker\delta\simeq\Ext^1_R(\Hom_R(M_1,N),R)$.
Regarding this diagram, we have the following three facts.
\roster
\item $\ker h_0\simeq \coker h_1$,
\item $\Tor^R_1(M,N^*)\simeq \ker h_1$, and
\item $\Ext^1_R(\Hom_R(M_1,N),R)\simeq \coker h_0$.
\endroster\medskip

Now assume that $\Ext^i_R(M,N)=0$ for all $i>0$, equivalently
$\Ext^1_R(M_i,N)=0$ for all $i\geq 0$. 
Constructing such a diagram as above for each of the short exact sequences
$0\to M_{i+1}\to F_i\to M_i\to 0$ and using the corresponding three 
facts as above we obtain exact sequences
$$
0\to\Tor_i^R(M,N^*)\to M_i\otimes_R N^* @>h_i>>\Hom_R(M_i,N)^*\to
\Tor^R_{i-1}(M,N^*)\to 0\tag{\bf 5.4}
$$
for $i\geq 2$. From Theorem 2.1 we know that $\Ext^i_R(M,N)=0$ for all
$i\gg 0$ if and only if $\Tor^R_i(M,N^*)=0$ for all $i\gg 0$.  Hence

\medskip

\proclaim{5.5 Proposition} Let $R$ be a local Gorenstein ring, and
let $M$ and $N$ be finitely generated $R$-modules with $N$ maximal
Cohen-Macaulay.  Then
$\Ext^i(M,N)=0$ for all $i\gg 0$ implies
the natural maps $M_i\otimes_RN^*\to\Hom_R(M_i,N)^*$ are isomorphisms
for all $i\gg 0$.\Qed
\endproclaim

Theorem 5.9 below contains a similar result.

\medskip

\subhead Ext and Stable Hom \endsubhead

Recall that the \it stable Hom \rm, $\underline{\Hom}_R(M,N)$,
is the cokernel of the natural map $M^*\otimes_R N \to \Hom_R(M,N)$.
Equivalently, it is the quotient of 
$\Hom_R(M,N)$ by maps $f:M\to N$ which factor
through a free module.  Stable Homs offer a convenient way of interpreting
the vanishing of all higher $\Ext^i_R(M,N)$: from (5.1) (and
the exact sequence involving Ext preceding it) we see that
$$
\Ext^i_R(M,N)\simeq\underline{\Hom}_R(M_i,N) \tag \bf 5.6 \rm
$$
for $i\geq 1$. Hence we may record the following as a corollary of 5.2.

\medskip

\proclaim{5.7 Corollary} Let $R$ be a local Gorenstein ring, and let
$M$ and $N$ be finitely generated $R$-modules with $M$ maximal 
Cohen-Macaulay.  Then $\Ext^i_R(M,N)=0$ for all $i \gg 0$ if and only if
for all $i \gg 0$ every map $M_i \to N$ factors through a free module.
\endproclaim

\medskip

The next Proposition allows us to shift among the stable Homs with ease,
which often can clarify basic vanishing results concerning Ext. 
(\it cf. \rm Remark 2.2 and (1.2).)
 
\proclaim{5.8 Proposition} Let $M$ and $N$ be finitely generated
maximal Cohen-Macaulay modules over the local Gorenstein ring $R$.  Then
we have
\roster
\item $\underline\Hom_R(M,N)\simeq\underline\Hom_R(N^*,M^*).$
\item $\underline\Hom_R(M,N)\simeq\underline\Hom_R(M_t,N_t)$ 
for all $t\in \Bbb Z$.
\endroster
\endproclaim
 
\demo{Proof} (1).  The isomorphism is induced by the obvious
mapping $\Hom_R(\text{ },R):\Hom_R(M,N) \to \Hom_R(N^*,M^*)$.
The fact that the induced map is an isomorphism is straightforward
(since $M$ and $N$ are reflexive) provided it is well-defined.  
But this is clear since if $
f \in\Hom_R(M,N)$ factors through a free module $F$, then
$\Hom_R(f,R)$ factors through $F^*$.
 
(2). It is enough to prove (2) in the case $t = 1$.
Given a map $f\in\Hom_R(M,N)$ we get a map
$f_1\in\Hom_R(M_1,N_1)$ by completing the diagram
$$
\CD
0 @>>> M_1 @>>> F @>\epsilon>> M @>>> 0 \\
@.     @VVf_1V  @VVf_0V        @VVfV  @.\\
0 @>>> N_1 @>>> G @>\epsilon'>> N @>>> 0.
\endCD\tag{5.8.1}
$$
Define $\Phi:\underline\Hom_R(M,N) \to\underline\Hom_R(M_1,N_1)$
by $\Phi(\bar f)=\bar f_1$.
 
We first show that $\Phi(\bar f)$ is determined independent of the
choice of chain map $\{f_0,f_1\}$.  Suppose that
$g_0:F \to G$ and $g_1:M_1\to N_1$ are two other maps making
the diagram (5.8.1) commute and such that $\epsilon'g_0=f\epsilon$.
Then we have the standard homotopy $h:F \to N_1$ such that
$f_1-g_1=h(M_1\hookrightarrow F)$.  That is, $f_1-g_1$ factors through
a free
module, so that $\bar f_1=\bar g_1$ in $\underline\Hom_R(M_1,N_1)$.
 
Next we show $\Phi$ is well-defined.
Suppose that $f\in\Hom_R(M,N)$ factors through a free
module $H$. Then $f_1=0$ completes the diagram
$$
\CD
0 @>>> M_1 @>>> F @>\epsilon>> M @>>> 0 \\
@.     @VVV  @VVV        @VVV  @.\\
{} @. 0 @>>> H @>>> H @>>> 0\\
@.     @VVV  @VVV        @VVV  @.\\
0 @>>> N_1 @>>> G @>\epsilon'>> N @>>> 0.
\endCD
$$
Hence, $\Phi(\bar f)=0$, as desired.
 
In order to show that $\Phi$ is an isomorphism we exhibit its inverse.
Let $g$ be in $\Hom_R(M_1,N_1)$.  We dualize and complete the diagram
$$
\CD
0 @>>> M^* @>>> F^* @>>> M_1^* @>>> 0 \\
@.     @AA(g^*)_1A  @AA(g^*)_0A        @AAg^*A  @.\\
0 @>>> N^* @>>> G^* @>>> N_1^* @>>> 0.
\endCD
$$
Define $\Psi:\underline\Hom_R(M_1,N_1)\to\underline\Hom_R(M,N)$
by $\Psi(\bar g)=\overline{(g^*)_1^*}$.  It's not hard to see that
$\Phi$ and $\Psi$ are inverses of one another (since $M$ and $N$ 
are reflexive).
\Qed
\enddemo
 
\medskip

\subhead Vanishing Ext and Cohen-Macaulayness \endsubhead

We end with a theorem which shows the relationship between the vanishing
of $d$ consecutive Ext modules and the Cohen-Macaulayness of a certain
tensor product. 
 
\proclaim{5.9 Theorem} Let $R$ be a $d$-dimensional local
Gorenstein ring, and let
$M$ and $N$ be maximal Cohen-Macaulay modules.
Consider the following two conditions.
\roster
\item $M^*\otimes_RN$ is maximal Cohen-Macaulay,
\item $\Ext^1_R(N,M)=\cdots=\Ext^d_R(N,M)=0$.
\endroster
Then (2) implies (1). If we assume that 
$\Ext^1_R(N,M),\dots,\Ext^d_R(N,M)$ have finite
length, then (1) implies (2).
Furthermore, if (1) holds then $\Hom_R(N,M)$ is maximal
Cohen-Macaulay and $M^*\otimes_RN \simeq
\Hom_R(N,M)^*$.  If (1) holds and $R$ is also integrally closed,  
then $M^*\otimes_RN \simeq\Hom_R(M,N)$.
\endproclaim
 
\demo{Proof} We first prove that if (1) holds then
the $R$-module $\Hom_R(N,M)$ is maximal Cohen-Macaulay
and $M^*\otimes_R N \simeq \Hom_R(N,M)^*$.  Note that as both
$M$ and $N$ are reflexive, the natural map
$\Hom_R(N,M)\to\Hom_R(M^*,N^*)$ is an isomorphism.
Since $M^*\otimes_RN$ is maximal Cohen-Macaulay so is
$(M^*\otimes_RN)^*$, and we have
$$
\alignat 2
(M^*\otimes_RN)^* &=\Hom_R(M^*\otimes_RN,R)&&{}\\
{}&\simeq\Hom_R(M^*,N^*)&&\qquad\text{by Hom-tensor adjointness}\\
{}&\simeq\Hom_R(N,M)&&{}\tag 5.9.1
\endalignat
$$
Hence $\Hom_R(N,M)$ is maximal Cohen-Macaulay and
$M^*\otimes_RN\simeq(M^*\otimes_RN)^{**}\simeq\Hom_R(N,M)^*$.  

We prove (1)$\Longrightarrow$(2) 
under the assumption that
$\Ext^1_R(N,M),\dots,\Ext^d_R(N,M)$ have finite
length. We induce on $d$.  
The case in which $d=0$ is vacuous.
 
$d = 1$. Since $\Ext^1_R(N,M)$ has finite length and $M$ and $N$ 
are maximal Cohen-Macaulay, we can choose a non-zerodivisor
$x\in R$ such that $x$ is a non-zerodivisor on both $M$ and $N$
and $x\Ext^1_R(N,M)=0$.  For any $R$-module $X$ we let $\overline X$
denote $X/xX$, and we let ${}^*$ indicate Hom into either $R$ 
or $\overline R$ depending on the module in question.  We have
the short exact sequence $0\to M @>x>> M \to \overline M\to 0$, which
yields the exact sequence
$$
0\to \Hom_R(N,M) @>x>> \Hom_R(N,M) \to \Hom_{\overline R}(\overline N,
\overline M)\to \Ext_R^1(N,M)\to 0.
$$
Hence $\length(\Hom_{\overline R}(\overline N,\overline M))=
\length(\overline{\Hom_R(N,M)})+\length(\Ext^1_R(N,M))$.
Note that for any maximal Cohen-Macaulay $R$-module $X$, 
$\overline{X^*}\simeq\overline X^*$. We have
$$
\alignat 2
\length(\Hom_{\overline R}(\overline N,\overline M)) &=\length((\overline M^*
\otimes_{\overline R} \overline N)^*)&&{}\qquad\text{by (5.9.1)}\\
{}&=\length(\overline M^* \otimes_{\overline R} \overline N)&&
\qquad\text{since $\overline R$ is 0-dimensional}\\
{}&=\length(\overline{M^*}\otimes_{\overline R}\overline N)&&
\qquad\text{since $\overline M^*\simeq\overline{M^*}$}\\
{}&=\length(\overline{M^*\otimes_R N})&&{}\\
{}&=\length(\overline{\Hom_R(N,M)^*})&&
\qquad\text{by (5.9.1)}\\
{}&=\length(\overline{\Hom_R(N,M)}^*)&&\\
{}&=\length(\overline{\Hom_R(N,M)})&&\qquad\text{since $\overline R$
is 0-dimensional}.
\endalignat
$$
Therefore $\length(\overline{\Hom_R(N,M)})=\length(\overline{\Hom_R(N,M)})
+\length(\Ext^1_R(N,M))$, \linebreak and so
$\Ext^1_R(N,M)=0$.
 
$d> 1$.  Choose a parameter $x\in\cap_{i=1}^d\ann_R\Ext^i_R(N,M)$.
We have $\overline M^*\otimes_{\overline R} \overline N \simeq
\overline{M^*\otimes_RN}$ is maximal Cohen-Macaulay. The short exact
sequence $0\to M @>x>> M \to \overline M\to 0$ yields  exact
sequences
$$
0\to\Ext^i_R(N,M)\to\Ext^i_R(N,\overline M)\to\Ext^{i+1}_R
(N,M)\to 0\tag 5.9.2
$$
for $i=1,\dots,d-1$.  Hence $\Ext^i_R(N,\overline M)$ have finite
length for $1\leq i\leq d-1$.  Since $\Ext^i_R(N,\overline M)
\simeq\Ext^i_{\overline R}(\overline N,\overline M)$ for all $i$ \cite{R}, it 
follows by induction that
$\Ext^1_{\overline R}(\overline N, \overline M)=
\cdots=\Ext^{d-1}_{\overline R}(\overline N,\overline M)=0$.  
Now the exact sequences (5.9.2), for $i=1,\dots,d-1$, and the fact that 
$\Ext^i_R(N,\overline M)
\simeq\Ext^i_{\overline R}(\overline N,\overline M)$ for all $i\geq 1$ gives 
$\Ext^i_R(N,M)=0$ for $i=1,\dots,d$.

Assume $(2)$. Note then that $\Hom_R(N,M)$ is maximal Cohen-Macaulay:
let 
$\F: \cdots \to F_1\to F_0\to N\to 0$ be an
$R$-free resolution of $N$. Applying $\Hom_R(\text{ },M)$ and
using our hypothesis we get the exact sequence
$$
0\to\Hom_R(N,M)\to\Hom_R(F_0,M)\to\cdots\to\Hom_R(F_{d+1},M).
$$
Each $\Hom_R(F_i,M)$ is maximal Cohen-Macaulay since $M$ is.
By counting depths along this exact sequence we get the
desired conclusion.

For $(2)\Longrightarrow (1)$ we again induce on $d$. 
The case $d= 0$ is trivial.
 
$d=1$.  Let $x$ be a non-zerodivisor on both $M$ and $N$.
The short exact sequence $0\to M @>x>> M\to\overline M\to 0$
and our hypothesis yield the short exact sequence
$$
0\to\Hom_R(N,M)@>x>>\Hom_R(N,M)\to\Hom_{\overline R}(\overline N,\overline M)
\to 0.
$$
Therefore 
$$
\overline{\Hom_R(N,M)}\simeq
\Hom_{\overline R}(\overline N,\overline M)\tag 5.9.3
$$
Consider the natural map $M^*\otimes_RN @>h>>\Hom_R(N,M)^*$.
We have a  commutative diagram
$$
\CD
\overline{M^*\otimes_RN} @>\overline h>>\overline{\Hom_R(N,M)^*}\\
@V\simeq VV @V\simeq VV\\
\overline M^*\otimes_{\overline R}\overline N @>\simeq>>
\Hom_{\overline R}(\overline N,
\overline M)^*
\endCD
$$ 
where the right vertical arrow comes from (5.9.3) and the
bottom arrow is the isomorphism of (5.9.1).
Thus $h$ is an isomorphism modulo $x$.  By Nakayama's lemma,
$h$ must be onto.  Now we have a short exact sequence
$$
0\to K\to M^*\otimes_R N @>h>>\Hom_R(N,M)^*\to 0.
$$
The fact that $\Hom_R(N,M)$ is maximal Cohen-Macaulay
implies that $\overline K=0$, and therefore $K=0$.
Thus $h$ is an isomorphism and $M^*\otimes_R N$ is 
maximal Cohen-Macaulay.

$d>1$.  For $x$ a non-zerodivisor on $M$ and $N$, the hypothesis
yields $\Ext^1_{\overline R}(\overline N,\overline M)=
\cdots =\Ext^{d-1}_{\overline R}(\overline N,\overline M)=0$,
so by induction $\overline M^*\otimes_{\overline R}\overline N$
is maximal Cohen-Macaulay, which means so is $M^*\otimes_R N$.

Finally suppose that $R$ is integrally closed. There is always a natural
map from $\Hom_R(M,N)\ra \Hom_R(N,M)^*$ obtained by composition, and
this map is an isomorphism if either $N$ or $M$ is free. Since $M$ and
$N$ are maximal Cohen-Macaulay modules and $R_P$ is regular if
the height of $P$ is at most one, it follows that this natural map
is an isomorphism in codimension one. It is a standard result that a 
map between reflexive modules which is an isomorphism in codimension one
must itself be an isomorphism.  Hence, as both $\Hom_R(M,N)$ and
$\Hom_R(N,M)^*$ are reflexive, the natural map map
$\Hom_R(M,N)\ra \Hom_R(N,M)^*$ is an isomorphism. The stated 
isomorphism of $M^*\otimes_RN$ with $\Hom_R(M,N)$ follows from the 
above paragraph. \Qed
\enddemo
\medskip
\head 6. Questions \endhead
\bigskip

This work leaves quite a few questions unresolved. We list a few for
further study.  Perhaps the most intriguing is
\medskip
1. Are all local Gorenstein rings AB rings?
\medskip
\noindent Some other interesting questions are:
\medskip
2. Let $R$ and $S$ be AB rings which are essentially of
finite type over a field $k$. Is $R\otimes_kS$ locally an AB ring?
\medskip
3. Are localizations of AB rings AB rings?
\medskip
4. Are AB rings Ext-bounded?


\enddocument